\input amsppt.sty
\magnification=1200  
\vsize=8.9truein
\hsize=6.5truein
\input epsf
   \def\norm#1{\left\Vert #1 \right\Vert}
   \def\rlhook#1{\lhook{\mathrel{\mkern-9mu}}
              \hbox to#1{\rightarrowfill}}

   \hskip-0.4truecm

\def\g{\Gamma}

      \font\sc=cmcsc10
\topmatter
\noindent Printing date : 12.8.96
\bigskip

\title
 Estimates for simple random walks \\ on fundamental groups of surfaces 
\endtitle

\author
Laurent Bartholdi, Serge Cantat, \\
Tullio Ceccherini Silberstein and Pierre de la Harpe 
\endauthor

\address
Section de Math\'ematiques, Universit\'e de Gen\`eve, C.P. 240,
CH-1211 Gen\`eve 24 (Suisse).\newline
E-mail: laurent.bartholdi\@epfl.ch \qquad pierre.delaharpe\@math.unige.ch
\endaddress

\address
Ecole Normale Sup\'erieure de Lyon, 46 All\'ee d'Italie,
69364 Lyon  Cedex 07 (France).\newline
   E-mail: Serge.Cantat\@ens.ens-lyon.fr
\endaddress

\address
Dipartimento di Matematica Pura e Applicata, Universit\`a degli Studi dell'
Aquila, \newline
Via Vetoio I-67100 L'Aquila (Italy).\qquad\qquad
   E-mail: tceccher\@mat.uniroma1.it  
\endaddress


\thanks
The third author acknowledges support from the 
\lq\lq Fonds National Suisse de la Recherche Scientifique\rq\rq. 
\endthanks 

\subjclass{ 60 J 15, 20 F 32}\endsubjclass
\abstract
   Numerical estimates are given for the spectral radius of
simple random walks on  Cayley graphs. Emphasis is on the
case of the fundamental group of a closed surface, for the usual system of
generators. 
\endabstract

\rightheadtext{Random walks on surface groups} 
\leftheadtext{Bartholdi - Cantat - Checcherini Silberstein - de la Harpe}

\endtopmatter

\document

\bigskip
\head
   Introduction
\endhead
\bigskip

Let $X$ be a connected graph, with vertex set $X^0.$ We denote by $k_x$ 
the number of neighbours of a vertex $x \in X^0.$ The {\it Markov operator} 
$M_X$ of $X$ is defined on functions on $X^0$ by
$$
   \left( M_Xf \right)(x) \ = \ \frac{1}{k_x} \sum_{y \sim x} f(y)
    \qquad f : X^0 \to \Bbb C \qquad x \in X^0
$$
where the summation is taken over all neighbours $y$ of $x$ 
(we assume that $1 \le k_x < \infty$ for all $x \in X^0$). \par

   If $X$ is a {\it regular graph,} i.e. if $k_x = k$ is independent on 
$x \in X^0,$ this operator
induces a  bounded self-adjoint operator on the Hilbert space $\ell ^2 (X^0),$
again  denoted by $M_X.$ The {\it spectral radius}
$
   \mu (X)
$
of the graph $X$ is the norm of this bounded operator. It is also a measure of
the asymptotic probability for a path of length $n$ 
in $X$ to be closed, and has several other interesting interpretations 
(see e.g. \cite{Woe}).
This carries over to the case of a not necessarily regular graph,
but the definition of the appropriate Hilbert space is slightly more
complicated (see again \cite{Woe}, Section 4.B). \par

   Let $\Gamma$ be a group generated by a finite set $S$ which is symmetric
$(s \in S \Longleftrightarrow s^{-1} \in S)$ and which does not contain the unit
element $1 \in \Gamma.$ Denote by $\text{Cay}(\g,S)$ the {\it Cayley graph} with
vertex set $X^0 = \Gamma$ and, for $x,y \in \Gamma,$ with $\{x,y\}$ an edge if
$x^{-1}y \in S.$ We denote by
$$
      \mu (\g , S)
$$ 
the spectral radius of the graph $\text{Cay}(\g,S).$ 
\par

   Let us recall two important results due to Kesten \cite{Ke1}, \cite{Ke2}. The
first one is the relation
$$
   \frac{ 2\sqrt{k-1}}{k} \ \le \ \mu (\g,S) \ \le \ 1
$$
with equality on the right if and only if $\g$ is amenable ($k$ is the number
of generators in $S$). For the second one let
us assume (for simplicity) that $\g$ does not have any element of order $2,$ so
that $k = 2h$ for some integer $h \ge 1;$ assume also (again for simplicity) that
$h \ge 2.$ Then one has the equality 
$$ 
    \frac{ \sqrt{2h-1}}{h} \ = \  \frac{ 2\sqrt{k-1}}{k} \ = \ \mu (\g,S) 
$$
if and only if $\g$ is a free group on a set $S_+ = \{s_1 ,..., s_h\}$ such that
$S = S_+ \coprod S_+^{-1}$ (where $\coprod$ indicates a disjoint union). \par

   There are few examples of exact {\it computations} of $\mu(\g,S)$ for non
amenable groups. Most of those we are aware of  are for groups which contain
free subgroups of finite index, even if there are a few known cases beyond
these \lq\lq almost free\rq\rq \ groups (see e.g. \cite{Car, Theorem 2} 
and \cite{CaM}). One  direction for further progress is to find good  
{\it estimates} for new classes of examples. \par

   As a test case, we consider here the fundamental group of an orientable 
closed surface of genus $g \ge 2,$ namely the group $\g_g$ given by the 
presentation
$$
   \g_g \ = \ \left\langle a_1 , b_1 , \hdots , a_g , b_g \ \bigg\vert \ 
    \prod_{j=1}^g a_jb_ja_j^{-1}b_j^{-1} \ = \ 1 \ \right\rangle
$$
and the generating set
$$
   S_g \ = \ \Big\{ \ a_1 \ , \ a_1^{-1} \ ,\  b_1 \ ,\  b_1^{-1} \ ,\ \hdots \ , \
        a_g \ , \ a_g^{-1} \ , \ b_g \ , \ b_g^{-1} \ \Big\}
$$
with $k = 4g$ elements; the resulting Cayley graph is denoted by $X_g.$ 
\par

Setting $\mu_g = \mu(X_g) = \mu(\g_g , S_g),$ one has
$$
   \frac{ \sqrt{4g-1}}{2g} \ < \ \mu_g \ < \ 1
$$
by Kesten's estimates recalled above. In particular
$$
   0.6614 \ \approx \ \frac{\sqrt{7}}{4} \ < \ \mu_2 \ < \ 1
$$
when $g = 2.$  As $\Gamma_g$ has $2g$ generators and as
$X_g$ has cycles of length $4g,$ the previous estimate may be improved to 
$$
    \frac{ \sqrt{4g-1}}{2g} \ + \ \frac{4 - 2 \sqrt 3}{(4g+2)(4g)^{4g+2}}
    \ \le \ \mu_g \ < \ 1 
$$
(see Formula (4.15) in \cite{Kes}), which gives for $g = 2$ an improvement 
of order $5 \times 10^{-11}.$  There is a better result due to Paschke, 
for which the improvement  is about $1.75 \ \times \ 10^{-4}$ \cite{Pas}. 
\bigskip
  
   In Section 1 below, we expose a very simple method based on an observation
of O. Gabber to show  that
$$
   \mu_g \ \le \ \frac{ \sqrt{2g-1}}{g} \qquad \text{and in particular} \qquad
   \mu_2 \ \le \ \frac{ \sqrt3 }{2} \ \approx \ 0.8660 .
$$
Section 2 records a computation with Poisson kernels; though it is in our view
the most interesting part of the present work, its numerical outcome so far is
limited to the inequality 
$$
       \mu_2 \ \le \ 0.7675
$$ 
and to similar inequalities for other
small values of $g.$ Section 3 uses embedding of trees in graphs to
improve the results of Section 1;  more precisely one has 
$$
   \mu_g \ \le \ \frac{ \sqrt{4g-2} }{2g} \ + \ \frac{1}{4g} 
   \qquad \text{and in particular} \qquad
   \mu_2 \ \le \ \frac{ \sqrt 6 }{4}  + \frac{1}{8} \ \approx \ 0.7373 .
$$
(One can extend much of Sections 1 and 3 to $C'(1/6)$ small cancellation groups
and to one relator groups.)
It follows from Section 3 and from Kesten's result that
$$
\mu_g \ = \ g^{-1/2} + O(g^{-1})
$$ 
for large $g.$
\bigskip


Our numerical results for $g \le 10$ are summarized in the following table.

$$\matrix
\text{genus}	&&& \text{Kesten}      	 && \text{Section 1}        && \text{Section 2}	&& 
                \text{Section 2}    && \text{Section 3}        \\
             &&&                      &&                         &&                  &&
                                    &&                         \\
	         g  &&& \frac{\sqrt{4g-1}}{2g}	&& \frac{\sqrt{2g-1}}{g}	&& 	 \nu  && 
               1 - \alpha	          && \frac{\sqrt{4g-2}}{2g} + \frac1{4g}  \\
             &&&                      &&                         &&                  &&
                                    &&                         \\
        --- 	&&& ---	                && ---	                     && ---	&& 
                ---	               && ---                                  \\
             &&&                      &&                         &&                  &&
                                    &&                         \\
2	 &&& .6614	&& .8660	&& .2990	&& .7675	&& .7373 \\
3	 &&& .5529	&& .7453	&& .2944	&& .6588	&& .6104 \\
4	 &&& .4841	&& .6615	&& .2932	&& .5872	&& .5303 \\
5	 &&& .4359	&& .6000	&& .2926	&& .5352	&& .4742 \\
             &&&                      &&                         &&                  &&
                            &&                                               \\
6 	&&& .3997	&& .5529	&& .2920	&& .4953	&& .4325 \\
7	 &&& .3712	&& .5153	&& .2916	&& .4633	&& .3999 \\
8 	&&& .3480	&& .4841	&& .2912	&& .4369	&& .3736 \\
9	 &&& .3287	&& .4581	&& .2908	&& .4147	&& .3518 \\
10	&&& .3123	&& .4359	&& .2905	&& .3956	&& .3332 \\
\endmatrix
$$
\medskip
   For example, for $g = 3,$ one has the lower bound $\mu_3 \ge 0.5529$ (Kesten)
and the upper bounds
$$
\alignat 2
\mu_3 \ &\le \ \frac{\sqrt 5}{3} \ \approx \ 0.7453 
          &&\qquad\text{(method of Section 1)}\\
\mu_3 \ &\le \ 0.6588 
          &&\qquad\text{(method of Section 2 with $\nu = 0.2944$)}\\
\mu_3 \ &\le \ \frac{\sqrt{10}}{6} + \frac{1}{12} \ \approx \ 0.6104  \ 
          &&\qquad\text{(method of Section 3).}
\endalignat
$$

\bigskip

   After completion of this work, the method of Section 1 has been improved by A. Zuk
\cite{Zuk}, who has shown in particular that
$$
   \mu_g \ < \ \frac{1}{\sqrt g}
$$
for all $g \ge 2,$ and again by T. Nagnibeda \cite{Nag}, who has shown in
particular that $$
   \mu_2 \ \le \ 0.6629.
$$
\bigskip

   We are grateful to Marc Burger, Bill Paschke, Rostislav Grigorchuk, Alain Valette 
and Wolfgang Woess for useful comments. 

\bigskip\bigskip
\head
 1.  Upper bounds from discrete $1$-forms
\endhead
\bigskip

   Let $X$ be a graph with vertex set $X^0$ and with edge set $X^1.$ Denote by 
$\Bbb X ^1$ the set of {\it oriented} edges of $X$ (if $X$ is finite, then
$\vert \Bbb X ^1 \vert = 2 \vert X ^1 \vert$). For each $e \in \Bbb X ^1$ we
denote by $\overline e$ the oriented edge obtained from $e$ by reversing the
orientation. A {\it $1$-form} on $X$ with values in some group $G$ is an
application $\omega : \Bbb X ^1 \to G$ such that 
$\omega(\overline e) = \omega(e)^{-1}$ for all $e \in \Bbb X ^1.$ 
We denote by $\Bbb R^*_+$ the multiplicative group $]0,\infty[.$ \par
   The following proposition is due to O. Gabber. It can be found in \cite{CdV}
(with the proof below) and its corollary in \cite{ChV} (with a different proof). 

\proclaim{Proposition 1} Let $X$ be a regular graph of degree $k.$ Suppose there
exists a $1$-form $\omega : \Bbb X ^1 \to \Bbb R^*_+$ and a constant $c > 0$ such
that
$$
   \frac{1}{k}\sum_{e \in \Bbb X ^1 , e_+ = x} \omega(e) \ \le \ c
$$
for all $x \in X^0.$ Then 
$$
      \mu(X) \ \le \ c .
$$
\endproclaim

\noindent   (The summation in the proposition holds over all oriented edges $e$      
heading to the vertex $x.$)

   \proclaim{Corollary 1} One has
$$
   \mu_g \ \le \ \frac{ \sqrt{2g-1}}{g}
$$
for all $g \ge 2.$ In particular
$$
  \mu_2 \ \le \ \frac{ \sqrt3 }{2} \ \approx \ 0.8660 .
$$
\endproclaim

\demo{Proof of Corollary 1} As the only relation in the chosen presentation of
$\g_g$ has even length, any edge $e$ in the Cayley graph $X_g$ of $(\g_g,S_g)$
joins two vertices $e_+ , e_-$ at different distances from the vertex $1.$
Let $d(x,y)$ denote the combinatorial distance in a graph between two vertices
$x,y,$ and write $\ell (x)$ for $d(1,x).$
For a number $b \ge 1$ (to be made precise below), one may thus define a $1$-form 
on $X_g$ by
$$
  \omega(e) \ = \ \left\{ \aligned 
     \  b^{-1} \qquad &\text{if} \qquad \ell (e_+) < \ell (e_-) \\
     \  b \ \ \qquad      &\text{if} \qquad \ell (e_+) > \ell (e_-) .
   \endaligned \right.
$$
Say that a vertex $x$ in $X_g$ is {\it of type $t$} if the set
$$
   \Big\{ \ y \in X_g \ \vert \ d(y,x) = 1 \ \ \text{and} \ \ 
     \ell (y) = \ell (x) - 1 \ \Big\}
$$
is of cardinality $t.$
For example $x$ is of type 1 if $0 < \ell (x) < 2g,$ and $x$ is of type $2$ if $x$
is at distance $2g$ from $1$ on a $4g$-gon containing $1.$ It follows from the
definition that $1$ is the only vertex of type $0.$ \par

   It is a fact that any other vertex is either of type $1$ or of type $2.$
This is well known and goes back to M. Dehn (or Poincar\'e ?); it is for example 
a straightforward consequence of Lemma 2.2 in \cite{Ser}. Compare with \cite{Can} 
and \cite{Wag}; note  however that a vertex is type $1$ [respectively type $2$] 
in our sense if and only if its Cannon type is in $\{1 , \hdots , 2g-1\}$   
[resp. is $2g$]. For convenience to the reader, we give a proof of the fact
we use in {\it Appendix A} below. \par

   One has
$$
   \sum_{e \in \Bbb X ^1 , e_+ = x} \omega(e) \ = \ 
   \left\{ \aligned
      \ 4gb^{-1} \ \ \ \ \   &\text{if} \ \  x = 1 
                       \ \ \text{(type 0)} \\
      \ (4g-1)b^{-1}+b \ \ \  &\text{if} \ \ x \ \ \text{is of type} \ \ 1 \\
      \ (4g-2)b^{-1}+2b \ \ &\text{if} \ \ x \ \ \text{is of type} \ \ 2
   \endaligned \right.
$$
and Proposition 1 applies with 
$$
 c \ = \ \frac{ (4g-2)b^{-1}+2b }{k}.
$$ 
To minimize $c,$ one sets $b = \sqrt{2g-1},$ so that
$$
  c \ = \ \frac{ 4\sqrt{2g-1}}{4g} .
$$
$\square$
\enddemo

\demo{Proof of Proposition 1} Let $f \in \ell^2(X^0).$
Choose $e \in \Bbb X ^1;$ set $x = e_+$ and $y = e_-.$ From
$$
   \left( \sqrt{\omega(e)} \vert f(x) \vert \ - \ 
          \frac{1}{\sqrt{\omega(e)}} \vert f(y) \vert \right)^2 \ \ge \ 0
$$
one has
$$
   2 \vert f(x) \vert \vert f(y) \vert \ \le \
   \omega(e) \vert f(x) \vert ^2 \ + \ \omega(\overline e) \vert f(y) \vert ^2 .
$$
Summing over $e \in \Bbb X ^1$ one obtains
$$
\aligned
   2 \sum_{x \in X^0} \vert f(x) \vert 
   \sum_{e \in \Bbb X ^1 , e_+ = x} \vert f(e_-) \vert \ &\le \ 
\\ \qquad
   \sum_{x \in X^0} \vert f(x) \vert ^2  
           &\sum_{e \in \Bbb X ^1 , e_+ = x} \omega(e) \ + \ 
   \sum_{y \in X^0} \vert f(y) \vert ^2  
           \sum_{e \in \Bbb X ^1 , \overline{e}_+ = y} \omega(\overline e)
\endaligned
$$
and
$$
 2 k \ \big\vert \left\langle f \mid M_X f \right\rangle \big\vert \ = \ 
 2 k \left\vert \sum_{x \in X^0} \overline{f(x)} \left(M_Xf\right)(x) \right\vert
 \ \le \ 2kc\norm{f}^2 .
$$
As this holds for all $f \in \ell^2(X^0),$ and as the operator $M_X$ on
$\ell^2(X^0)$ is self-adjoint, one has $\norm{M_X} \le c$ and the
conclusion follows.
$\square$ 
\enddemo

  {\bf Generalization.} {\it Let  
$\Gamma = \left\langle S_+ \vert R \right\rangle$
be a  group presentation satisfying a small cancellation hypothesis $C'(1/6).$
If $h \doteq \vert S_+ \vert \ge 2$ and if $S = S_+ \cup \left(S_+\right)^{-1}$, 
one has
$$
    \mu(\Gamma,S) \ \le \ \frac{ 2\sqrt{h-1}Ê}{h} \ .
$$
}

\demo{Proof} One has $\vert S \vert = 2h$ because small cancellation groups
cannot have elements of order $2$ (see e.g. Section V:4 in \cite{LyS}).
Types being defined as in the proof of Corollary 1, it is known that any vertex
distinct from the identity in the Cayley graph of $(\Gamma,S)$ is either of type
$1$ or of type $2$ (lemme 4.19 in \cite{Cha}). Defining a 1-form $\omega$ on this
Cayley graph by 
$$
   \omega(e) \ = \ 
   \left\{ \aligned
       b^{-1} \ \ \ \ \    &\text{if} \ \ \ \ell(e_+) \ < \ \ell(e_-) \\
   \ \   1    \ \ \ \ \ \ \ \ \    &\text{if} \ \ \ \ell(e_+) \ = \ \ell(e_-) \\
   \ \   b    \ \ \ \ \  \ \ \ \   &\text{if} \ \ \ \ell(e_+) \ > \ \ell(e_-) 
   \endaligned \right.
$$
one may apply verbatim the argument of Corollary 1.
$\square$ 
\enddemo

\bigskip\bigskip
\head
 2.  Upper bounds from Poisson kernels
\endhead
\bigskip

   Let again $X = \text{Cay}(\g , S)$ be as in the Introduction and let $M_X$ be
the corresponding Markov operator. The {\it combinatorial Laplacian} of $X$ is
defined to be
$$
   \Delta _X \ = \ 1 - M_X .
$$
Let $\alpha \in \Bbb R;$ a function $f : \g \longrightarrow [0,\infty[$ is 
said to be {\it $\alpha$-superharmonic} if $f \ne 0$ and if 
$\Delta_X f \ge \alpha f.$ (If there exists such a function $f$, one has
$f \ge \Delta_Xf \ge \alpha f$ and consequently $\alpha \le 1.$ One may also
show that $f(\gamma) > 0$ for all $\gamma \in \g.$) The function is said to
be {\it $\alpha$-harmonic} if moreover $\Delta_X f = \alpha f.$ 

\proclaim{Proposition 2} Let $\alpha \in \Bbb R.$ The following are equivalent.
$$
\aligned
   (i) \ \ & \ \alpha \le 1 - \mu(X) = \inf \Big\{ \text{ spectrum of }\ \Delta_X\
               \text{ on the Hilbert space } \ell^2(\g) \ \Big\}. \\
   (ii) \ \ & \text{There exists a function } f : \g \ \longrightarrow [0,\infty[
                \text{ which is $\alpha$-superharmonic} . \\
   (iii) \ \ & \text{There exists a function } f : \g \ \longrightarrow [0,\infty[
                \text{ which is $\alpha$-harmonic} . 
\endaligned
$$
\endproclaim

   There is one proof in terms of graphs in \cite{DoK, Proposition 1.5}. 
But there are earlier proofs in the literature on irreducible stationary
discrete Markov chains; the equivalence of (i) and (ii) is standard; 
the equivalence with (iii) is more delicate (see \cite{Har} and
\cite{Pru}). 

\proclaim{Corollary 2} One has $\mu_2 \ \le\  0.784.$ \par
   More generally, upper estimates for $\mu_g$ and small $g$ 's are given
by the table in the introduction. \endproclaim

  We begin the proof of Corollary 2 with the following lemma.

\proclaim{Lemma 1} Let $g$ be an integer, $g \ge 2.$ Set
$$
   D_g \ = \ 2 \operatorname{\arg\,cosh} \left( \cot \frac{\pi}{4g} \right) 
   \tag1
$$
For $\phi \in [0, 2\pi [$, set
$$
   b(\rho , \phi) \ = \ 
     \frac{ 1 }{ \cosh \rho \ - \ \sinh \rho \ \cos \phi }  \tag2
$$
for all $\rho > 0$ and
$$
  F_g(\nu , \phi) \ = \ \frac{1}{4g}  \sum_{j=0}^{4g-1} \ 
   \left\{ b \left( D_g \ , \ \phi+j\frac{2\pi}{4g}\right) \right\}^{\nu}
   \tag3
$$
for all $\nu \in \Bbb R.$ Then
$$
   \mu_g \ \le \  \max_{0 \le \phi < 2\pi}F_g(\nu , \phi) 
$$
for all $\nu \in \Bbb R$.
\endproclaim

   \demo{Proof}\enddemo

   \demo{First step : definition of a function $f_{\nu}$} \smallskip  

   Let $H^2$ be the hyperbolic plane.\par
   
   There is a free discrete isometric action of $\g_g$ on $H^2$ and a point 
$z_0 \in H^2$ such that the Dirichlet cells of the orbit $\g_g z_0$ constitute a
tesselation of $H^2$ by regular $4g$-gons with all inner angles equal to 
$\frac{\pi}{2g}.$ There is consequently an embedding of the graph
$X_g = \text{Cay}(\Gamma_g,S_g)$ in $H^2,$ vertices of the graph corresponding
to points of the orbit $\Gamma_g z_0$ and edges of the graph to pairs
of adjacent Dirichlet cells. Trigonometric computations for a hyperbolic triangle
with angles $\pi/2 , \pi/4g , \pi/4g$ show that $D_g$ in $(1)$ is the distance
between the centres of two adjacent Dirichlet cells.\par
   
   Let $\omega_0 \in \partial H ^2$ be a point at infinity. Let $P : H^2 \to
]0,\infty[$  be the function given by the value at $\omega_0$ of the Poisson
kernel. For computations we choose 
$$
   H^2 \ = \ \Big\{ \ z \in \Bbb C \ \big\vert \ Im(z) > 0 \ \Big\} \   
   \ \ \  \text{and}  \ \ \  \omega_0 \ = \ \infty i
    \ \ \  \text{so that}  \ \ \  P(x+iy) \ = \ y . \tag4
$$
Let $\Delta_H$ be the hyperbolic Laplacian on $H^2.$ One has 
$$
   	\Delta_H P^{\nu} \ = \ - \nu ( \nu - 1 ) P^{\nu}
$$
for all $\nu \in \Bbb R.$ (We have chosen a {\it positive} Laplacian $\Delta_H.$
This implies that the spectrum of the corresponding self-adjoint
operator on the Hilbert space  
$L^2 \left (H^2 , y^{-2}dxdy \right)$ is $[\frac{1}{4},\infty[.$
The equality $\Delta_H P^{\nu}  =  - \nu ( \nu - 1 ) P^{\nu}$ shows that there
exists $\alpha$-harmonic functions for $\Delta_H$ for all
$\alpha \le \frac{1}{4},$ in accordance with an analogue  for $\Delta_H$ of the
previous proposition. Much more on this in \cite{Sul}.) \par

   We define
$$
   f_{\nu} \ : \ \g_g \ \longrightarrow \ ]0,\infty[
$$
by $f_{\nu}(\gamma) = P^{\nu}(\gamma z_0).$ For $\gamma \in \g,$ let
$z_{\gamma , j}$ ($0 \le j \le 4g-1$) denote the centers of the Dirichlet cells
adjacent to the Dirichlet cell centered at $\gamma z_0.$ One has
$$
   \left( \Delta_X f_{\nu} \right) (\gamma) \ = \ 
   P^{\nu}(\gamma z_0 ) \ - \ 
   \frac{1}{4g} \sum_{j=0}^{4g-1} P^{\nu} \left(z_{\gamma,j}\right)
$$
for each $\gamma \in \g.$ The strategy of the proof is to find some
$\alpha \in \Bbb R$ such that $\Delta_X f_{\nu} \ge \alpha f_{\nu},$
and to deduce from the previous proposition that $\mu_g \le 1 - \alpha.$ 
\enddemo

\demo{Second step : lower estimate for $\Delta_Xf_{\nu}$} \smallskip

    For $z \in H^2 \ , \ \rho > 0$ and $\phi \in [0,2\pi[,$ let 
$z(\rho , \phi) \in H^2$ be the point at hyperbolic distance $\rho$ from $z$ for
which the oriented angle between the geodesic ray $\overarrow{z_0,\omega_0}$ 
and the geodesic segment $\overarrow{z_0,z(\rho , \phi)}$ is $\phi.$ Set
$$
   c_g(\nu , \rho , \phi , z) \ = \ \frac{ 
       P^{\nu}(z) \ - \  \frac{1}{4g} \sum_{j=0}^{4g-1} 
         P^{\nu} \bigg( z \left(\rho , \phi + j\frac{2\pi}{4g}
                  \right) \bigg)      }
      {   P^{\nu}(z)  }  .  \tag5
$$
Observe that there is one well-defined value 
$\phi_{\gamma} \in [0, \frac{2\pi}{4g}[$ such that
$$
  \left( \Delta_X f_{\nu} \right) (\gamma) \ = \ 
   c_g(\nu , D_g , \phi_{\gamma} , \gamma z_0 ) \ f_{\nu}(\gamma)
$$
for each $\gamma \in \Gamma.$ But computing the angles $\phi_{\gamma}$ is a  
difficult task, and we rather look for an estimate of the right-hand side
in the inequality
$$
  \Delta_X f_{\nu}  \ \ge \ \left( \min_{ 0 \le \phi < 2\pi \atop z \in H^2} 
   c_g(\nu , D_g , \phi , z) \right) \ f_{\nu}.  
$$  
Now $(5)$ shows that $c_g(\nu , \rho , \phi , z)$ depends neither on the real
part  of $z,$ because $P(x+iy) = y$ for all $x \in \Bbb R,$ nor on the
imaginary  part  of $z,$ because $P^{\nu}(\lambda z) = \lambda ^{\nu}
P^{\nu}(z)$ for all  $\lambda > 0.$ Thus one has
$$
   \Delta_X f_{\nu} \ \ge \ 
   \left( \min_{0 \le \phi < 2\pi} c_g(\nu , D_g , \phi , z_0) \right) \ f_{\nu}.
$$
Choosing moreover $z_0 = i,$ one has
$$
   P(z_0) \ = \ 1
$$
and
$$
   c_g(\nu , D_g , \phi , z_0) \ = \ 1 \ - \ \frac{1}{4g} \sum_{j=0}^{4g-1}
    \left\{ \Im \bigg( z_0 \Big( D_g , \phi + j\frac{2\pi}{4g} \Big) \bigg) 
                           \right\}^{\nu}
$$
by $(5).$
\enddemo

\demo{Third step : computation of $\Im \left( z_0(\rho , \phi) \right)$}
\smallskip

   Let $\Cal C$ be a hyperbolic circle of hyperbolic radius $\rho$ centered at the
point $z_0 = i$ of the Poincar\'e half-plane. The Cartesian coordinates $(a,b)$
of a point on $\Cal C$ satisfy
$$
      a^2 \ + \ (b - \cosh \rho)^2 \ = \ \left( \sinh \rho \right)^2 .  \tag6
$$
For each $\phi \in ]-\pi , \pi[,$ let $\Cal C_{\phi}$ be the hyperbolic geodesic
through $z_0$ defining at this point an angle $\phi$ with the vertical axis. The
Cartesian coordinates of a point on $\Cal C_{\phi}$ satisfy
$$
    \left( a - \frac{1}{\tan \phi} \right)^2 + b^2 \ = \ 1 + 
      \frac{1}{\tan^2 \phi } \ . \tag7
$$

Let us compute the second coordinates of the two points of 
$\Cal C \cap \Cal C_{\phi}$ (see Figure 1).
Subtracting $(7)$ from $(6),$ one finds
$$
   \frac{a}{\tan \phi} \ -\  b \ \cosh \rho \ = \ -1
$$
and inserting this in $(7)$ one obtains
$$
   \Big( \cosh^2 \rho \ \tan^2 \phi + 1 \Big)  b^2 \ - \ 
   2 \Big(\cosh \rho \ ( \tan^2 \phi + 1 ) \Big) b \ + \ 1 + \tan^2 \phi \ = \ 0 .
$$
Straightforward manipulations show that 
$$
  \Big(\cosh \rho \ (\tan^2 \phi \ + \ 1 ) \Big) ^2 \ - \ 
   \Big( \cosh^2 \rho \ \tan^2 \phi + 1 \Big)\Big( 1 + \tan ^2 \phi \Big) \ = \ 
   \left( \frac{ \sinh \rho}{\cos \phi} \right) ^2 
$$
and consequently that
$$
\aligned
        b \ = \ \frac{ 
         \cosh \rho \ (\tan^2 \phi + 1) \ \pm  \ \frac{\sinh \rho}{ cos \phi }  }
        { \cosh ^2 \rho \ \tan^2 \phi \ + \ 1  } 
      \ &= \  \frac{ \cosh \rho \ \pm  \ \sinh \rho \ cos \phi }
              { \cosh ^2 \rho \ \sin^2 \phi \ + \ \cos^2 \phi } \\
      \ &= \ \frac{ 1 }{ \cosh \rho \ \mp \ \sinh \rho \ \cos \phi }. 
\endaligned \tag8
$$
Thus one has
$$
     \Im\Big( z_0(\rho , \phi) \Big) \ = \ 
   \frac{ 1 }{ \cosh \rho \ - \ \sinh \rho \ \cos \phi } \ = \ b(\rho , \phi) 
$$
where the last equality is $(2).$ (The other sign in (8) would give 
$b(\rho , \phi + \pi)$.) \enddemo

\smallskip
\centerline{\epsfbox{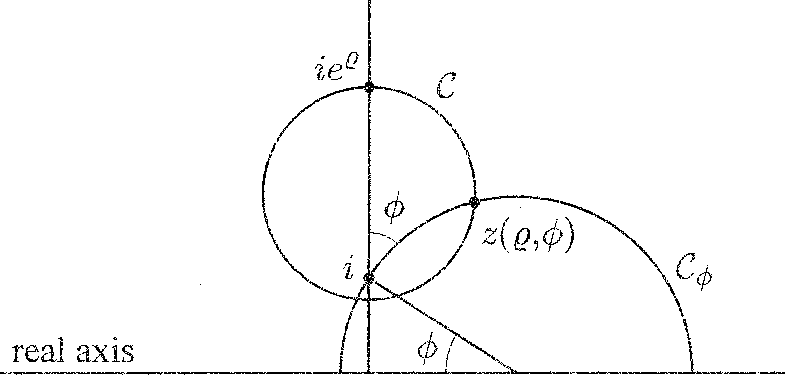}}
$$ \text{ Figure 1. } $$
\smallskip

\demo{Fourth step : coda}\smallskip

   The previous computations show that one has
$$
   \Delta_X f_{\nu} \ \ge \ \alpha f_{\nu}
$$
for
$$
   \alpha = \min_{0 \le \phi < 2\pi} \Big\{ 1 \ - \ F_g(\nu , \phi) \Big\} 
$$
where $F_g$ is defined in $(3).$
As $\mu_g \le 1 - \alpha$ by Proposition 2, this ends the proof.
$\square$
\enddemo

\medskip

   At this point, the problem is to compute
$\inf_{\nu} \max_{\phi} F_g(\nu , \phi).$
One could use just here a computer system such as {\sc Maple} 
and obtain a table of numerical results. However we rather
adopt the following program.\par

   A first step consists of a lemma of calculus showing that,
for any $\nu \in [0,1],$ the function $\phi \mapsto F_g(\nu,\phi)$
reaches its maximum at $\phi = 0.$ (This at least for $g \le 27;$ we have 
not found a reasonably short proof working for all $g.$) This is
stated below, and proved in the {\it Appendix B} at the end of our paper. \par

   Only in a second step we use a computer, first to find an
efficient value of $\nu$ (which turns out to be near $0.3$ for all $g$)
and then to compute $F_g(\nu,0)$ for this $\nu,$ so that one has a
numerical estimate
$$
   \mu_g \ \le \ F_g(\nu , 0)
$$
for the spectral radius of 
$\mu_g = \mu\left(\text{Cay}(\Gamma_g,S_g)\right).$

\bigskip

   For $g$ and $\nu$ fixed, the function $\phi \mapsto 4gF_g(\nu,\phi)$
is a sum of a function 
$$
        \beta : \phi \mapsto \left(\cosh (D_g) - \sinh (D_g) \cos \phi \right)^{-\nu}
$$
and of $4g-1$ translates of $\beta.$
It is straightforward to check that $\dot \beta (0) = 0$ and 
$\ddot \beta (0) < 0$, so that $\beta$ has a local maximum at the origin. 
The purpose of Lemma 2 (which is proved in  Appendix B) is to show that this local
maximum is  strong enough for $\phi \mapsto F_g(\nu,\phi)$ to have an absolute
maximum at the origin.

\proclaim{Lemma 2} For $2 \le g \le 27$ and $0 \le \nu \le 1$ one has
$$
   \max_{0 \le \phi \le 2\pi} F_g(\nu , \phi) \ = \ F_g(\nu , 0) . 
$$
Thus, for these $g$'s,
$$
   \mu_g \ \le \ F_g(\nu , 0) \ 
$$
for all $\nu \in [0,1],$ by Lemma 1.
\endproclaim

  \demo{End of proof of Corollary 2} Thanks to the previous lemma,
we may consider the function
$$
   \nu \ \ \longmapsto \ \ 
  F_g(\nu,0) \ = \ \frac1{4g}\sum_{j=0}^{4g-1}\beta\left(j\frac{2\pi}{4g}\right),
$$
and compute its minimum over $0\le\nu\le 1$, yielding an upper bound for
$\mu_g$. The computer algebra program {\sc Maple} was used here, giving for
$g \le 10$ the  values of the table in the Introduction.
$\square$
\enddemo

\bigskip\bigskip
\head
 3.  Upper bounds from regular subtrees
\endhead
\bigskip

   Let $X$ be a regular graph of degree $k,$ as in Section 1. Assume that there
is a subgraph $Y$ of $X$ which is spanning (namely which contains all vertices 
of $X$) and which is regular of degree $l$ for some $l \in \{2 , \hdots , k-1 \}$
(we assume $k \ge 3$). The Markov operators $M_X$ and $M_Y$ act on the same
space  $\ell ^2(X^0) = \ell ^2(Y^0).$ One has
$$
\aligned
   \left( M_X f \right) (x) \ &= \ 
    \frac{1}{k} \left\{ \sum_{e \in \Bbb Y ^1 , e_+ = x} f(e_-) \ + \ 
         \sum_{e \in \Bbb X ^1 \setminus \Bbb Y^1, e_+ = x} f(e_-) \right\} \\ 
       &= \ \frac{l}{k} \left( M_Yf \right) (x) \ + \ \frac{1}{k} 
          \sum_{e \in \Bbb X ^1 \setminus \Bbb Y^1, e_+ = x} f(e_-)
\endaligned
$$
so that
$$
   \norm{M_X} \ \le \ \frac{l}{k} \norm{M_Y} \ + \ \frac{k-l}{k} .
$$
In case $Y$ is a disjoint union of regular trees, $\norm{M_Y}$ is explicitely
known from Kesten's computations and one has the following.

   \proclaim{Proposition 3} Let $X$ be a regular graph of degree $k \ge 3$ and
let $Y$ be a spanning subgraph of $X$ which is a disjoint union of regular
trees of degree $l,$ for some 
\newline 
$l \in \{2 , \hdots , k-1 \}.$ Then
$$
   \frac{ 2 \sqrt{k-1} }{k} \ \le \ \norm{M_X} \ \le \ 
   \frac{ 2 \sqrt{l-1}Ê}{k} \ + \ \frac{k-l}{k} .
$$
\endproclaim

   \proclaim{Lemma 3} The graph $X_g$ contains a spanning subgraph $Y_g$ which
is a disjoint union of regular trees of degree $4g - 1.$
\endproclaim

\demo{Proof} Recall from Section 1 that $\ell (x)$ denotes the combinatorial
distance in $X_g$ between a vertex $x$ and the base point $1,$ 
and  from {\it Appendix A} that vertices in $X_g$ are
shared amongst three {\it types} numbered $0, 1$ and $2.$ Recall also that
\smallskip

   (a) two vertices of type 2 are at distance at least $3$ from each other,
\smallskip

   (b) any vertex $x$ of type $1$ has a {\it convenient} neighbour 
$y \in X^0_g$   such that  \par
 \qquad  $\ell (y) = \ell (x) + 1,$ \par
  \qquad $y$ is of type $1,$ \par
  \qquad all neighbours of $y$ in $X_g$ are of type $1$ \par
[indeed $x$ has at least $4g-2$ of these neighbours].

\smallskip \noindent The construction goes in two steps. \smallskip

   {\it First step.} Let $Z_g$ be the spanning subgraph of $X_g$ obtained from
$X_g$ by erasing, for each vertex $x$ of type $2,$ one edge connecting $x$ to a
neighbour $y$ of $x$ such that $\ell (y) = \ell (x) -1.$ (This edge is chosen
arbitrarily from $2$ candidates.) By (a) above, any vertex of type $1$ has
degree $4g-1$ or $4g$ in $Z_g$ and any vertex of type $2$ has degree $4g-1$ in
$Z_g.$
\smallskip

   {\it Second step.} For each $k \ge -1,$ define inductively a graph
$Y_g^{(k)}$ as follows. First, set $Y_g^{(-1)} = Z_g.$
Then, if $k \ge 0,$ let $Y_g^{(k)}$ be a spanning subgraph of $X_g$ obtained from
$Y_g^{(k-1)}$ by erasing, for each vertex $x$ with $\vert x \vert = k$ which
is of  degree $4g$ in $Y_g^{(k-1)},$  one edge connecting $x$ to one of its
convenient neighbours. 
(This edge is chosen arbitrarily from at least $4g-2$ candidates.)  
By (b) above, any vertex with $\vert x \vert \le k$ in $Y_g^{(k)}$ is of degree
$4g-1.$ \par
   Observe that, for all $l \ge k,$ the graphs $Y_g^{(k)}$ and $Y_g^{(l)}$
coincide \lq\lq in the ball defined by $\vert x \vert \le k$\rq\rq . Thus one may
set $Y_g = Y_g^{(\infty)};$ any vertex in $Y_g$ is of degree $4g-1.$
\smallskip

   Let us check that $Y_g$ does not contain any circuit. For this, we will show
that $Z_g$ has no circuit. \par

   Observe that two neighbours in $Z_g$ are never at the same distance from
$1$ (because this is already so in $X_g,$ a consequence of the relation
defining the group $\Gamma_g$ being of {\it even} length). If  there were a
circuit in $Z_g,$ it would contain a vertex $x$ at maximum distance, say $n$,
from $1,$ and this $x$ would have two neighbours at distance $n-1;$ in
particular $x$ would be of type $2;$ this is ruled out by the first step above.
\par

   Thus $Y_g$ is indeed a spanning forest of degree $4g-1$ in $X_g.$ \smallskip

   Though this fact is not needed for what follows, let us observe that $Y_g$ has
infinitely many connected components. Indeed, choose a vertex $x$ of type $1$
and a convenient neighbour $y$ of $x$ such that the edge connecting $x$ to $y$
has been erased in the second step above; then any neighbour $z$ of $y$ in
$Y_g$ is such that $\ell (z) = \ell(y) + 1.$ Choose similarly a vertex
$x' \ne x$ and a convenient neighbour $y',$ with the same properties as $x$ and
$y.$ Then $y$ and $y'$ are not in the same component of $Y_g,$ because any path
from $y$ to $y'$ in $Y_g$ should have a maximum strictly between $y$ and $y',$
and this is ruled out by the first step above.\par
   There are infinitely many such $x$'s, because from (a) there are infinitely many
vertices of type $1$ and degree $4g$ in $Z_g.$ 
$\square$ \enddemo

   {\it Remark.} In another terminology, Lemma 3 shows that the set of edges of
$X_g$ which  are not edges of $Y_g$ constitute a {\it perfect matching} of $X_g,$
also called a {\it $1$-factor.}

   \proclaim{Corollary 3}  One has
$$
   \mu_g \ \le \ \frac{ \sqrt{4g-2} }{2g} \ + \ \frac{1}{4g}
$$
for all $g \ge 2.$ In particular 
$$
   \mu_2 \ \le \ \frac{ \sqrt 6  }{4}  + \frac{1}{8} \ \approx \ 0.7373 .
$$
\endproclaim

\demo\nofrills{Proof:\usualspace} immediate from Proposition 3 and Lemma 3.
$\square$
\enddemo

   {\it Comparison with Corollary 1.} Computations in this section are more
efficient that computations of Section 1 (with discrete $1$-forms), because
$$
   \frac{ \sqrt{4g-2} }{2g} \ + \ \frac{1}{4g}   
                             \ < \ \frac{\sqrt{2g-1}}{g}
$$
for all $g \ge 2.$ But Computations of Section $1$ can be improved to beat the
present ones \cite{Nag} !

   \proclaim{Corollary 4} Let 
$\Gamma = \left\langle S_+ \big\vert R \right\rangle$
be a one-relator group, with $S_+ \subset \Gamma \setminus \{1\}$ of order 
$h \ge 2.$ Then
$$
  \frac{ \sqrt{2h-1}Ê}{h} \ < \ \mu(\Gamma,S) \ \le \
  \frac{ \sqrt{2h-3} + 1}{h} 
$$
for $S = S_+ \cup \left(S_+\right)^{-1}.$
\endproclaim

\demo{Proof} Let $T_+$ be a subset obtained from $S_+$ by erasing one letter
appearing in $R$ (we assume $R$ to be cyclically reduced). Then $T_+$ is free
by the Dehn-Magnus' Freiheitssatz (see e.g. \cite{ChM, Chapter II.5}).
Set $T = T_+ \cup \left(T_+\right)^{-1}.$ Let $Y$ be the spanning subgraph of the
Cayley graph $Cay(\Gamma,S)$ for which two vertices $x,y$ are connected by an
edge whenever $xy^{-1} \in T.$ As $T_+$ is free in $\Gamma,$ the graph $Y$ is a
disjoint union of regular trees of degree $2h-2.$ The corollary follows from
Proposition 3.
$\square$
\enddemo

\bigskip\bigskip
\head
   Appendix A : on planar graphs
\endhead
\bigskip

   Let $X$ be a connected graph embedded in the plane, edges of $X$ being 
piecewise smooth curves which are pairwise disjoint (but for common
vertices). If $X$ is infinite, we assume that the following {\it strong
planarity} condition holds: for any simple closed curve in $X,$ the
corresponding bounded region of the plane (via the Jordan curve theorem)
contains only {\it finitely many} vertices of $X.$ A {\it face} of $X$ is
the closure of a connected component of the complement of $X$ in the plane.
\par

   Let  $d(x,y)$ denote the combinatorial distance between two vertices 
$x,y \in X^0;$ let $x_0 \in X^0$ be a base point and set 
$\ell (x) = d(x_0,x).$ If $X$ is bipartite, two
neighbouring vertices $x,y \in X^0$ are necessarily such that
$\vert \ell (x) - \ell (y) \vert = 1.$ Recall that the {\it type} $t(x)$ of
a vertex $x \in X^0$ is {\it here} the number of neighbours $y$ of $x$ such that
$\ell (y) < \ell (x).$ Observe that, for $x \in X^0,$ one has $t(x) = 0$ if and
only if $x = x_0.$ \par

   \proclaim{Geometric proposition} Let $X$ be a strongly planar graph with
base point $x_0 \in X^0.$ Assume that $X$ is connected, bipartite, and
satisfies the following conditions: \smallskip

   (i) = large degree: each vertex $x \in X^0$ has $k_x \ge 4$ neighbours
in $X;$ \par
   (ii) = large faces: each face $F$ of $X$ contains $k_F \ge 4$ vertices
of $X;$ \par
   (iii) = no-sink-vertex: each vertex $x \in X^0$ has at least one
neighbour \par
  \qquad $y \in X^0$ such that $\ell (y) = \ell (x) + 1.$
\smallskip
\noindent Then $t(x) \le 2$ for all $x \in X^0.$ \smallskip

   Assume moreover that each face $F$ of $X$ contains $k_F \ge 8$
vertices of $X.$ Then \smallskip

 \noindent  (a) for two vertices $x,y$ of type $t(x) = t(y) = 2,$ one
has    $d(x,y) \ge 3,$ \par
  \noindent (b) any vertex $x$ of type $1$ has a neighbour $y \in X^0$
such that $d(x_0,y) = d(x_0,x) + 1$ \par\qquad
   and such that all neighbours of $y$ are also of type $1.$
\endproclaim

\demo{Proof}
   We will make use of the following {\it maximum principle}: if $C$ is a
simple closed curve in $X$ enclosing a bounded open region $R$ of the
plane, then
$$
   \max_{x \in R \cap X^0} d(x_0,x) \ < \ 
      \max_{y \in C \cap X^0} d(x_0,y) .
$$
To show this, consider a point $x' \in R$ and a geodesic segment from
$x_0$ to $x'.$ By (iii), this can be extended to an arbitrarily long
geodesic segment starting at $x_0.$ By strong planarity, such an 
extension has to escape $R$ and does so crossing $C$ in some vertex $y'.$
One has clearly $d(x_0,x') < d(x_0,y'),$ and this proves the inequality
above. \medskip

   We will also make use of another standard fact: for two distinct faces $F$ and
$G,$ the intersection $F \cap G$ is either empty, or a vertex of the graph, or one 
edge of the graph. (To rule out the case of several edges, one may evaluate the
Euler characteristics of the closure of a  bounded component of  the complement of
$F \cup G.$) \medskip

   {\it Claim A. For each face $F$ of $X,$ the function
$$
   f_F \ : \  \left\{ \aligned
        F \cap X^0 \quad &\longrightarrow \qquad \Bbb N \\
         x     \qquad &\longmapsto \quad \ell(x)
  \endaligned \right.
$$
has a unique local minimum (say $m_F$) and a unique local maximum (say $M_F$). 
In other words, the function $f_F$ is unimodal.}\smallskip

   To prove the claim, it is enough to show that, for any $n \in \Bbb N,$ the
cardinal of the fiber $f_F^{-1}(n)$ is at most $2.$ \par

\smallskip
\centerline{\epsfbox{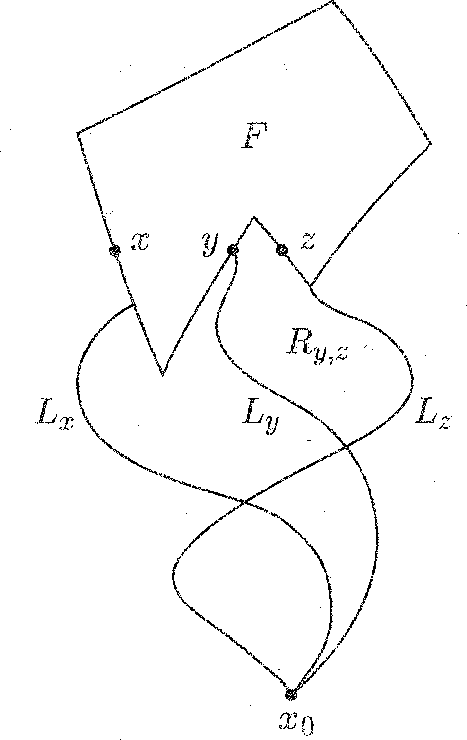}}
$$ \text{ Figure 2. } $$
\smallskip

  Suppose {\it ab absurdo} that this is not the case. Let
$x,y,z \in F \cap X^0$ be three distinct vertices such that
$f_F(x) = f_F(y) = f_F(z).$ Denote by $[x,y], [y,z], [z,x]$ the three sides of a
triangle with vertices $x,y,z$ contained in the boundary of $F.$ Choose geodesic
segments $L_x,L_y,L_z$ from $x_0$ to $x,y,z$ respectively. Then appropriate
subsegments of $[x,y],L_x,L_y$ constitute a simple closed curve $C_{x,y}$
defining a bounded open region $R_{x,y}$ of the plane; one has similarly curves
$\ C_{y,z}\ ,\ C_{z,x}\ $ and regions $\ R_{y,z}\ ,\ R_{z,x}\ .$ Let $R$ be the
interior of 
$\overline{R_{x,y}} \cup \overline{R_{y,z}} \cup \overline{R_{z,x}}.$ \ 
There is exactly one of the three points $x,y,z$ which is inside $R;$ upon
changing notations for $x,y,z,$ one may assume that $y \in R$ (as in Figure 2).

\par

   The geodesic segment $L_y$ can be extended indefinitely, by (iii).
Such an extension of $L_y$ has to escape $R$ through its boundary, and this is
impossible; thus  {\it Claim A} is proved. \smallskip

   It follows that  the two geodesic segments
in $F \cap X$ from $m_F$ to $M_F$ have the same number
$\ell(M_F) - \ell(m_F) -1$ of interior vertices - this number being
strictly positive by (ii). \medskip

   {\it Claim B. There is no vertex $x \in X^0$ with type $t(x) \ge 3.$}
\smallskip

   Indeed, suppose {\it ab  absurdo} that $X$ has vertices of type at
least $3$ and let $m$ be one of these for which the distance to $x_0$ is
minimum. Let $v_1 , \hdots , v_r , w_1 , \hdots , w_s$ be the neighbours
of $m,$ listed in such a way that
$$
\aligned
 \ell(v_i) \ &= \ \ell(m) - 1 \qquad 1 \le i \le r \qquad (r \ge 3),\\
 \ell(w_k) \ &= \ \ell(m) + 1 \qquad 1 \le k \le s \qquad (s \ge 1).
\endaligned
$$
For $i \in \{1 , \hdots , r\},$ choose a geodesic segment $L_i$ from
$x_0$ to $v_i.$ \par

\smallskip
\centerline{\epsfbox{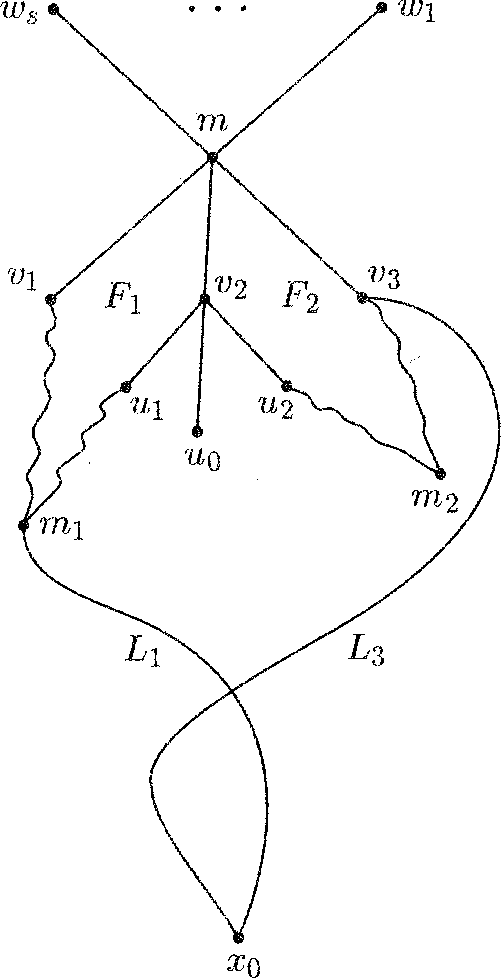}}
$$\text{Figure 3.}$$
\smallskip

   For $i,j \in \{1 , \hdots , r\}$ with $i \ne j,$ the segment
$[v_i , m , v_j]$ and appropriate subsegments of $L_i,L_j$ constitute a
simple closed curve $C_{i,j}$ defining a bounded open region $R_{i,j}$ of
the plane. By the maximum principle, $w_k \notin R_{i,j}$
for all $k \in \{1 , \hdots , s\}.$ Thus, upon renumbering the $v_i$ 's
and the $w_k$ 's, one may assume that 
$v_1 , \hdots , v_r , w_1 , \hdots , w_s$ are arranged in cyclic order
around the vertex $m.$ It follows that there is a face $F_1$ containing
$v_1,m,v_2,$ a face $F_2$ containing $v_2,m,v_3,$ and that $F_1,F_2$ are
adjacent along $[v_2,m]$ (see Figure 3). \par

   For $h \in \{1,2\},$ let $u_h$ denote the vertex of $F_h$ such that 
$d(u_h,v_2) = 1$ and $\ell(u_h) = \ell(v_2) -1;$ let also $m_h$ denote
the vertex of $F_h$ nearest to $x_0$ and choose a geodesic segment 
$\tilde L_h$ from $x_0$ to $m_h.$ (We have used {\it Claim A} here.) By (i), 
the vertex $v_2$ has a neighbour  $u_0 \in X^0 \setminus \{m,u_1,u_2\}.$ Using
again the maximum principle for a region enclosed by appropriate subsegments
of  $\tilde L_1 \cup [m_1,v_2]$ and $\tilde L_2 \cup [m_2,v_2],$ one checks 
that $\ell(u_0) = \ell(v_2) -1.$ It follows that $v_2$ is of type at
least $3$ (because it has neighbours $u_0,u_1,u_2$), in contradiction with
the choice of $m$ (because $\ell(v_2) < \ell(m)$); thus {\it Claim B} is
proved. \medskip

   {\it Proof of (a).}
Let $x,y \in X^0$ be such that $x \ne y$ and $t(x) = t(y) = 2.$
There is a face $F$ such that $x$ is the vertex of $F$ maximizing the
distance to the origin on $F \cap X^0,$ and a face $G$
associated similarly to $y.$
   The equality $d(x,y) = 1$ would contradict Claim B, as it is indicated
in Figure 4 (this uses only $k_H \ge 6$ for all faces $H$ of $X$). 

\smallskip
\centerline{\epsfbox{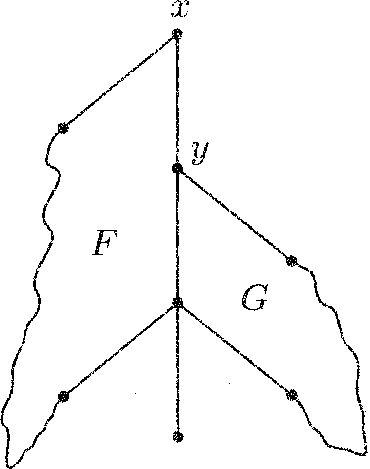}}
$$ \text{ Figure 4. } $$
\smallskip

The equality $d(x,y) = 2$ gives
rise to two type of configurations, each in contradiction with Claim B,
as it is indicated in Figure 5.

  {\it Proof of (b).}
Let $x \in X^0$ be a vertex of type $1.$ Let 
$v,w_1 , \hdots , w_s$ be the neighbours of $x,$ listed in cyclic order
around the vertex $x,$ with
$$
\aligned
   \ell(v) \   &= \ \ell(x) - 1 \\
   \ell(w_k) \ &= \ \ell(x) + 1  \qquad 1 \le k \le s \qquad (s \ge 3).
\endaligned
$$
We leave it to the reader to check the following facts : \par
   the vertices $w_1$ and $w_s$ are of types $1$ or $2$ (not both of type
$2$ by Claim B), \par
   the intermediate vertices $w_2 , \hdots , w_{k-1}$ are all of type
$1,$ \par
   any of these has all its neighbours of type $1.$ \smallskip

   This ends the proof of the proposition. $\square$

\enddemo

\smallskip
\centerline{\epsfbox{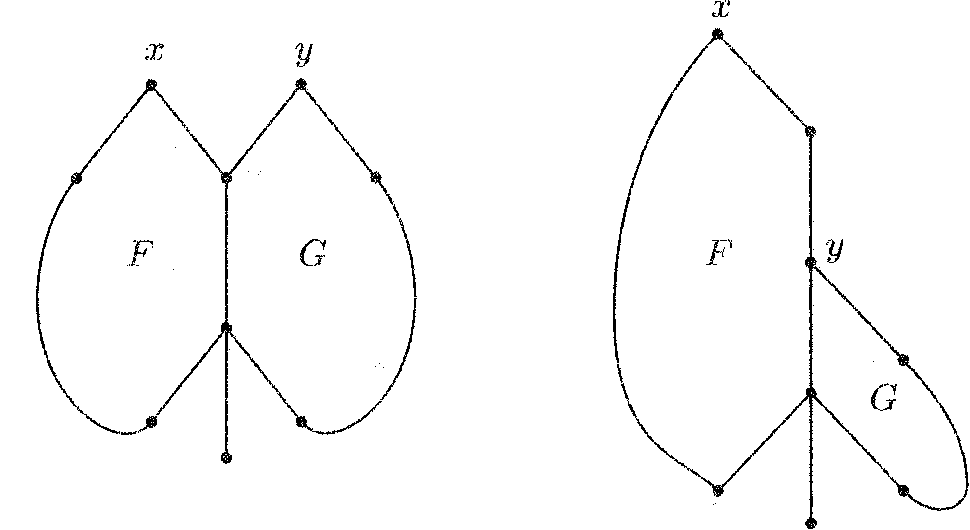}}
$$ \text{ Figure 5. } $$
\smallskip


\bigskip\bigskip
\head
Appendix B : proof of Lemma 2
\endhead
\bigskip

   \proclaim{Lemma 4} For $g \ge 2$, set
$$
\aligned
      C_g = \cosh (D_g) \qquad\qquad 
          &\delta_g \ = \ \arccos \left( \frac{S_g}{C_g} \right) 
           \ = \ \arccos \left( \tanh (D_g) \right) \\
     S_g = \sinh (D_g)  \qquad\qquad
        &\epsilon_g \ = \ \arccos \left( \frac{S_g}{C_g}
            - \frac{1}{S_gC_g} \right) \ .
\endaligned
$$
Then one has
$$
   0 \ < \ \delta_g \ < \ \epsilon_g \ < \ \frac{\pi}{4g} \ ,   \tag9
$$
and
$$
\alignat 3
 \frac{d}{d \phi} \ \frac {1}{\left(C_g - S_g \cos \phi \right) ^{\nu}}
       & \ \le \ 0 
 &
 \qquad &\text{for all} 
 &
 \qquad \phi \ &\in \ [0 , \pi]  \\
 \frac{d^2}{d \phi ^2} \ \frac {1}{\left(C_g - S_g \cos \phi \right) ^{\nu}}
       &\ \ge 0 \ 
 &
 \qquad &\text{for all} 
 &
 \qquad \phi \ &\in \ [\delta_g , \pi]  \\
 \frac{d^3}{d \phi ^3} \ \frac {1}{\left(C_g - S_g \cos \phi \right) ^{\nu}}
       & \ \le \ 0  
 &
 \qquad &\text{for all} 
 &
 \qquad \phi \ &\in \ [\epsilon_g , \pi].
\endalignat
$$
\endproclaim

\demo{Proof of Lemma 4}\enddemo

   \demo{First step : inequalities of (9) in Lemma 4}
Obviously $0 < \delta_g$, as $C_g$ and $S_g$ are both positive. Better,
$S_g > 1$ because $D_g > 1$; indeed $D_g$ is an increasing function of
$g$ (being the composite of two decreasing functions and an increasing one),
and $D_2 \approx 3.057 > 1$.
This allows us to write $S_g > S_g - 1/S_g > 0$; dividing
by $C_g$ and taking arccosines yields $\delta_g < \epsilon_g$.

Next $\epsilon_g < \pi/4g$. For this, as `$\cos$' is decreasing, we must show
that
$$
\frac{S_g}{C_g} - \frac1{S_gC_g} 
\   \overset{?}\to{>}  \
\cos\left(\frac\pi{4g}\right)                     \tag{10}
$$
holds without the $?$ sign.
We set $X=\cot^2(\pi/4g)$ and we express $C_g$, $S_g$, $\cos(\pi/4g)$ in
terms of $X;$ as 
$ \ \ C_g = \cosh(D_g) = 2\left(\cosh(\frac{D_g}{2})\right)^2 - 1 \ ,$
one has
$$
C_g = 2X-1 \qquad 
S_g = 2\sqrt{X(X-1)} \qquad 
\cos\left(\frac\pi{4g}\right)=\sqrt{\frac X{X+1}}
$$
whence (10) becomes
$$\frac{2\sqrt{X(X-1)}}{2X-1} - \frac1{2\sqrt{X(X-1)}(2X-1)} 
\   \overset{?}\to{>}  \ 
\sqrt{\frac X{X+1}}.
$$
Squaring,
$$
4X(X-1) - 2 + \frac1{4X(X-1)} 
\   \overset{?}\to{>}  \ 
\frac X{X+1}(2X-1)^2
$$
or, provided $X>1$,
$$
16X^4 - 44X^3 + 20X^2 + 9X + 1 
\   \overset{?}\to{>}  \ 
0.
$$
We rewrite this as
$$
16(X-2)^4 + 84(X-2)^3 + 140(X-2)^2 + 73(X-2) + 3 
\   \overset{?}\to{>}  \
0.
$$
This inequality is true for all $X>2$, as the left hand side is
a polynomial in $X-2$ with all coefficients positive.
It remains to check that $\cot^2(\pi/4g) > 2$ for all $g$; but this is clear
because $\cot^2(\pi/4g)$ is an increasing function of $g$ 
with value $3+2\sqrt2$ at $g=2$.
\enddemo

\demo{Second step: the function $\beta$} Set
$$
   \beta(\phi) \ = \ b(D_g,\phi)^{\nu} \ = \ 
   \frac{1}{ (C_g - S_g \cos \phi)^{\nu}} 
$$
so that 
$$
   F_g(\nu , \phi) \ = \ \frac{1}{4g} \sum_{j = 0}^{4g - 1} 
           \beta \left(\phi + j\frac{2\pi}{4g} \right) . \tag{11}
$$
The first derivative of $\beta$ is
$$
\dot \beta (\phi) \ = \ 
       \frac{-\nu S_g \sin \phi}{ (C_g - S_g \cos \phi)^{\nu + 1}} \tag 12
$$
so that $\dot \beta (\phi) \le 0$ for all $\phi \in [0,\pi].$
The second derivative of $\beta$ is
$$
\ddot \beta (\phi) \ = \ 
   \nu S_g \ \frac{ S_g - C_g \cos \phi + \nu S_g \sin^2 \phi }  
                    { (C_g - S_g \cos \phi)^{\nu + 2}}
\ \ge \  
   \nu S_g \ \frac{ S_g - C_g \cos \phi }  
                    { (C_g - S_g \cos \phi)^{\nu + 2}} \  \tag 13
$$
so that $\ddot \beta (\phi) \ge 0$ as soon as $\cos \phi \le S_g / C_g,$
namely as soon as $\phi \in [\delta_g , \pi].$
The third derivative of $\beta$ is
$$
\aligned
\dddot \beta (\phi) \ &= \ \nu S_g \sin \phi \ 
    \frac{ 1 \ - \ (3 \nu + 1)S_g \left(S_g - C_g \cos \phi \right) 
          \ - \ \nu ^2 S_g ^2 \sin ^2 \phi }
          { (C_g - S_g \cos \phi)^{\nu + 3}} \\ 
       &\le \   \nu S_g \sin \phi \ 
    \frac{ 1 \ - \ (3 \nu + 1)S_g \left(S_g - C_g \cos \phi \right) }
          { (C_g - S_g \cos \phi)^{\nu + 3}}
\endaligned
$$
so that $\dddot \beta (\phi) \le 0$ for $\phi \in [\epsilon_g , \pi].\square$
\enddemo

\bigskip

\demo{Proof of Lemma 2} Let $g \ge 2$ and $\nu \in [0,1]$ be fixed. 
As the function $\phi \mapsto F_g(\nu,\phi)$ is smooth, even
and periodic of period $\frac{\pi}{2g}$ it is enough to show that
$$
   F_g(\nu , \phi) \ \le \ F_g(\nu , 0)
$$
for all $\phi \in [0,\frac{\pi}{4g}].$ \medskip

   In the range $[\delta_g , \frac{\pi}{2g} - \delta_g],$
the functions 
$\phi \mapsto b\left( D_g , \phi + j\frac{2\pi}{4g}\right) ^{\nu}$
are convex for all $j \in \{0,1,\hdots,4g-1\}$
by Lemma 4. Their convex sum $\phi \mapsto F_g(\nu,\phi)$
is thus also convex, so that
$$
   F_g(\nu , \phi) \ \le \ F_g(\nu , \delta_g)
$$
for all $\phi \in [\delta_g , \frac{\pi}{4g}].$ \medskip

   We suppose now $\phi \in [0,\delta_g]$ and we want to show that
$\frac{d}{d\phi}F_g(\nu,\phi) \le 0.$  One has
$$
   \frac{d}{d\phi} \ 4gF_g(\nu , \phi) \ = \ \dot \beta (\phi) \ + \ 
       \sum_{j=1}^{4g-1} \dot \beta \left(\phi + j \frac{\pi}{2g} \right)
$$
by (11). As $\dot \beta$ is an odd function 
$\sum_{j=0}^{4g-1} \dot \beta \left( j \frac{\pi}{2g} \right) = 0;$ as 
$\dot \beta (0) = \dot\beta(\pi) = 0$ one has also
$$
   \frac{d}{d\phi} \ 4gF_g(\nu , \phi) \ = \ \dot \beta (\phi) \ + \ 
       \sum_{j=1}^{4g-1} \bigg( \dot \beta \left(\phi + j \frac{\pi}{2g} \right)
        \ - \   \dot \beta \left( j \frac{\pi}{2g} \right) \bigg)  . 
$$
By the theorem of Rolle,
$$
   \frac{d}{d\phi} \ 4gF_g(\nu , \phi) \ = \ \dot \beta (\phi) \ + \ 
       \sum_{j=1}^{4g-1} \phi \ddot \beta \left(\psi _j + j \frac{\pi}{2g} \right)
$$
for some $\psi _j \in [0,\phi].$ By the computation for $\dddot \beta$ in 
Lemma 4, one has
$\ddot \beta (\psi_j + j\frac{\pi}{2g}) \le \ddot \beta (\frac{\pi}{2g})$ 
and
$$
   \frac{d}{d\phi} \ 4gF_g(\nu , \phi) \ \le \ \dot \beta (\phi) \ + \ 
       (4g-1) \  \phi \ddot \beta \left( \frac{\pi}{2g} \right)  . 
$$
Using (12) and (13) one finds
$$
\multline
\frac{d}{d\phi} \ 4gF_g(\nu , \phi) \ \le \\
     -\nu S_g \frac{\sin \phi}{ (C_g - S_g \cos \phi)^{\nu + 1}} \ + \ 
     (4g-1) \nu S_g \phi \frac{ S_g - C_g \cos(\pi/2g) + \nu S_g \sin^2(\pi/2g) }
                  { \left( C_g - S_g \cos(\pi/2g) \right)^{\nu + 2}  }
\endmultline
$$
so all we have to check is
$$ 
     \frac{ \sin \phi / \phi} { (C_g - S_g \cos \phi)^{\nu + 1}}
     \ \ge \  (4g-1)   \frac{ S_g - C_g \cos(\pi/2g) + \nu S_g \sin^2(\pi/2g) }
                  { \left( C_g - S_g \cos(\pi/2g) \right)^{\nu + 2}  }
$$ 
for all $\phi \in [0,\delta_g].$

As $\cos\phi \ge \cos(\pi/2g)$, so $(C_g - S_g\cos\phi)^\nu \le
  (C_g - S_g\cos(\pi/2g))^\nu$, we may tighten the inequality to
$$
\frac{\sin\phi/\phi}{C_g - S_g\cos\phi} \ \ge \ 
  (4g-1)\frac{S_g - C_g\cos(\pi/2g)+\nu S_g\sin^2(\pi/2g)}{
    \left(C_g - S_g\cos(\pi/2g)\right)^2} \ \doteq \ R_g(\nu);
$$
as the right hand side is constant in $\phi$ 
while the left hand side decreases monotonically,
we let $\phi=\delta_g$. Finally we set $\nu=1$ to maximize the 
right hand side. Our goal is now to show
$$
\frac{\sin\delta_g/\delta_g}{C_g - S_g\cos\delta_g} \  \ge \ R_g(1).
$$
But, by definition of $\delta_g$ (see Lemma 4), one has
$C_g - S_g \cos \delta_g = 1/C_g$
and
$C_g \sin \delta_g = \sqrt{ C_g^2 - S_g^2} = 1,$
so that our goal reduces to showing
$$
   \frac{1}{\delta_g} \ \ge \ R_g(1).
$$
That this is true for $g \le 27$ can in turn be checked on a pocket calculator.
Thus when $g\le 27$ and $\nu\in[0,1]$ the function
$F_g(\nu,-)$ is monotonously decreasing
on $[0,\pi/4g]$; its maxima are at $0 + j\pi/2g$ and its minima at 
$\pi/4g + j\pi/2g$.
$\square$
\enddemo

\newpage


\Refs
\widestnumber\no{BCH2}



\ref \no Can \by J.W. Cannon \paper The growth of the closed surface groups and
compact hyperbolic Coxeter groups \jour Circulated typescript 
\yr University of Wisconsin, 1980 \endref

\ref \no Car \by  D.I. Cartwright \paper  Some examples of random walks
on free products of discrete groups \jour Ann. Mat. Pura Appl. \vol 151 \yr 1988
\pages 1--15 \endref

\ref \no CaM \by D.I. Cartwright and W. Mlotkowski \paper Harmonic analysis for
groups acting on triangle buildings \jour J. Austral. Math. Soc., Ser. A \vol 56
\yr 1994 \pages 345--383 \endref

\ref \no Cha \by C. Champetier \paper Propri\'et\'es statistiques des groupes de
pr\'esentation finie \jour Adv. in Math. \vol 116 \yr 1995 \pages 197--262 \endref

\ref \no ChM \by B. Chandler and W. Magnus \book The history of combinatorial
group theory  : a case study in the history of ideas \publ Springer \yr 1982
\endref

\ref \no ChV \by P.A. Cherix and A. Valette \paper On spectra of simple random 
walks on one-relator groups \jour Pacific J. Math. \vol 175 \yr 1996
\pages 417--438 \endref

\ref \no CdV \by Y. Colin de Verdi\`ere \paper Spectres de graphes  
\jour Cours Spécialisés, 4. Société Mathématique de France, Paris \yr
1998. \pages viii+114 pp \endref

\ref \no DoK \by  J. Dodziuk and L. Karp \paper Spectra and function theory for
combinatorial Laplacians \jour Contemp. Math. \vol 73 \yr 1988 \pages 25--40
\endref

\ref \no FP \by W.J. Floyd and S.P. Plotnick \paper Symmetries of planar growth
functions \jour Invent. Math. \vol 93 \yr 1988 \pages 501--543 \endref

\ref \no Har \by  T.E. Harris \paper Transient Markov chains with stationary 
measures \jour Proc. Amer. Math. Soc. \vol 8 \yr 1957 \pages 937--942 \endref

\ref \no Ke1 \by H. Kesten \paper Symmetric random walks on groups
\jour Trans. Amer. Math. Soc. \vol 22 \yr 1959 \pages 336--354 \endref

\ref \no Ke2 \by H. Kesten \paper Full Banach mean values on countable groups
\jour Math. Scand. \vol 7 \yr 1959 \pages 146--156 \endref

\ref \no LyS \by R.C. Lyndon and P.E. Schupp \book Combinatorial group theory
\publ Springer \yr 1977 \endref


\ref \no Nag \by T. Nagnibeda \paper  An upper bound for the spectral
 radius of a random walk on surface groups {\rm (Russian)} \jour
 Zap. Nauchn. Sem. S.-Peterburg. Otdel. Mat. Inst. Steklov. (POMI)
 \vol 240  \yr 1997 \pages Teor. Predst. Din. Sist. Komb. i Algoritm. Metody. 2, 154--165, 293--294; translation in  J. Math. Sci. (New York) {\bf 96}  (1999),  no. 5, 3542--3549\endref

\ref \no Pas \by  W.B. Paschke \paper Lower bound for the norm of a
vertex-transitive graph \jour Math. Zeit. \vol 213 \yr 1993 \pages 225--239
\endref


\ref \no Pru \by W.E. Pruitt \paper Eigenvalues of non-negative matrices \jour
Ann. Math. Stat. \vol 35 \yr 1964 \pages 1797--1800 \endref

\ref \no Ser \by C. Series \paper The infinite word problem and limit sets of
Fuchsian groups \jour Ergod. Th. \& Dynam. Sys. \vol 1 \yr 1981 \pages 337--360 
\endref

\ref \no Sul \by D. Sullivan \paper Related aspects of positivity in
Riemannian geometry \jour J. Diff. Geom. \vol 25 \yr 1987 \pages 327--351 \endref

\ref \no Wag \by P. Wagreich \paper The growth function of a discrete group 
\jour Springer Lecture Notes in Math. \vol 956 \yr 1982 \pages 125--144 \endref

\ref \no Woe \by W. Woess \paper Random walks on infinite graphs and groups --- a
survey on selected topics \jour Bull. London Math. Soc. \vol 26 \yr 1994 
\pages 1--60 \endref

\ref \no Zuk \by A. Zuk \paper A remark on the norm of a random walk on surface
groups \jour  Colloq. Math. \vol 72 \yr 1997 \pages 195--206\endref

\endRefs
\enddocument